\documentclass[11pt, leqno, a4paper]{article}
\usepackage{amscd, amsmath}
\usepackage{amssymb}
\usepackage{moreverb}
\def\noproof{{\unskip\nobreak\hfill\penalty50\hskip2em\hbox{}%
     \nobreak\hfill$\Box$\parfillskip=0pt%
     \finalhyphendemerits=0\par}}
\def\enddemo{\ifmmode\eqno\Box\else\noproof\vskip0.8truecm\fi}

\newtheorem{theo}{Theorem}

\newtheorem{remark}[theo]{Remark}
\newtheorem{remarks}[theo]{Remarks}
\newtheorem{lemma}[theo]{Lemma}
\newtheorem{proposition}[theo]{Proposition}

\newcommand{\lra}{\longrightarrow}

\DeclareMathOperator{\Spec}{Spec}

\DeclareMathOperator{\Res}{Res}

\DeclareMathOperator{\res}{res}
\DeclareMathOperator{\cor}{cor}
\DeclareMathOperator{\charac}{char}
\DeclareMathOperator{\id}{id}
\DeclareMathOperator{\Hom}{Hom}

\DeclareMathOperator{\Alb}{Alb}
\DeclareMathOperator{\CH}{CH}
\DeclareMathOperator{\Obj}{Obj}
\DeclareMathOperator{\Ker}{Ker}
\DeclareMathOperator{\Jac}{Jac}

\DeclareMathOperator{\Tr}{Tr}

\DeclareMathOperator{\Gal}{Gal}

\DeclareMathOperator{\Mod}{Mod}

\DeclareMathOperator{\ab}{\Mod_{\bZ}}
\DeclareMathOperator{\punkt}{\bf{\cdot}}
\DeclareMathOperator{\rank}{rank}

\DeclareMathOperator{\eff}{eff}
\DeclareMathOperator{\et}{et}
\DeclareMathOperator{\Nis}{Nis}
\DeclareMathOperator{\DM}{DM}
\DeclareMathOperator{\DMn}{\DM^{\eff,\mbox{--}}_{\Nis}}
\DeclareMathOperator{\DMe}{\DM^{\eff,\mbox{--}}_{\et}}

\newcommand{\Dnis}[1]{{D^{-}(Shv_{\Nis}(Cor_{#1}))}}
\newcommand{\Det}[1]{{D^{-}(Shv_{\et}(Cor_{#1}))}}

\newcommand{\bN}{{\mathbb N}}
\newcommand{\bQ}{{\mathbb Q}}

\newcommand{\bZ}{{\mathbb Z}}

\newcommand{\bG}{{\mathbb G}}

\newcommand{\cC}{{\cal C}}
\newcommand{\cO}{{\cal O}}
\newcommand{\wcO}{{\widehat{\cO}}}

\newcommand{\io}{{\iota}}

\newcommand{\noi}{\noindent}

\begin{document}

\title{A counterexample to generalizations of the Milnor-Bloch-Kato conjecture\footnote{This work was done during the second author stayed at Universit{\"a}t Bielefeld supported by SFB 701. He is grateful to the members there.}}
\author{By Michael Spie{\ss} and Takao Yamazaki}
\date{}
\maketitle

\begin{abstract}
We construct an example of a torus $T$ over a field $K$
for which the Galois symbol 
$K(K; T,T)/n K(K; T,T) \to H^2(K, T[n]\otimes T[n])$ 
is not injective for some $n$.
Here $K(K; T,T)$ is the Milnor $K$-group attached to $T$
introduced by Somekawa. We show also that the motive $M(T\times T)$ gives a counterexample to another generalization of the Milnor-Bloch-Kato conjecture (proposed by Beilinson).
\end{abstract}

\section{Introduction} 

Let $K$ be a field, $m$ a positive integer and $n$ an integer prime to the characteristic of $K$. The Milnor-Bloch-Kato conjecture asserts that the Galois symbol  
\begin{equation}
\label{eqn:galoisclass}
K^M_m(K)/n K^M_m(K) \quad \lra \quad H^m(K, \bZ/n\bZ(m))
\end{equation}
from Milnor $K$-groups to Galois cohomology is bijective. Recently, Rost and Voevodsky have announced a proof (special cases have been obtained earlier by Merkurjev-Suslin, Rost and Voevodsky). 

In \cite{somekawa}, Somekawa has introduced certain {\it generalized Milnor $K$-groups} $K(K; A_1, \ldots, A_m)$ attached to semi-abelian varieties $A_1, \ldots, A_m$. If $A_1 = \ldots = A_m = \bG_m$ is the one-dimensional split torus they agree with the usual $K^M_m(K)$. If $m=2$, $A_1 = \Jac_X$ and $A_2 = \Jac_Y$ are the Jacobians of smooth, projective and connected curves $X$ and $Y$ over $K$ having a $K$-rational point, then $K(K; A_1, A_2)$ is the kernel of the Albanese map $\CH_0(X\times Y)_{\deg = 0} \to \Alb_{X\times Y}(K)$. 

Somekawa has defined a Galois symbol 
\begin{equation}
\label{eqn:generalgaloisclass}
K(K; A_1, \ldots, A_m)/n K(K; A_1, \ldots, A_m) \lra H^m(K,A_1[n]\otimes\ldots\otimes A_m[n])
\end{equation}
and conjectured that it is always injective. In this note we present a counterexample (see section \ref{section:somekawa}). Let us describe it briefly. Let $L/K$ be a cyclic extension of degree $n$ and $\sigma$ a generator of the Galois group $\Gal(L/K)$. Let $T$ be the kernel of the norm map $\Res_{L/K} \bG_m \to \bG_m$. We show that the norm $K(L; T,T) \to K(K; T,T)$ induces an isomorphism  $K_2(L; T,T)/(1-\sigma) \to K_2(K; T,T)$. On the other hand, the corresponding map of Galois cohomology groups $H^2(L,T[n]\otimes T[n])/(1-\sigma) \to H^2(K,T[n]\otimes T[n])$ is neither injective nor surjective (for a suitable choice of $L/K$). Note that, since $T$ is split over $L$, the Galois symbol $K_2(L; T,T)/n \to  H^2(L,T[n]\otimes T[n])$ is bijective. Consequently, $K_2(K; T,T)/n \to  H^2(K,T[n]\otimes T[n])$ is in general not injective.

In the section \ref{section:beilinson} we show that the motive $M(T\times T)$ gives a counterexample to another generalization  of the Milnor-Bloch-Kato conjecture (proposed by Beilinson). 

We would like to thank Bruno Kahn for his comments on the first version of this note. In particular he pointed out to us that our counterexample to Somekawa's conjecture should also provide a counterexample to Beilinson's conjecture.

\section{Counterexample to Somekawa's conjecture}
\label{section:somekawa}

\paragraph{Algebraic groups as Mackey-functors} Let $K$ be a field. For a finite field extension $L/K$ and commutative algebraic groups $G$ over $K$ and $H$ over $L$ we denote by $G_L$ the base change of $G$ to $L$ and by $\Res_{L/K} H$ the Weil restriction of $H$. The functor $G \mapsto G_L$ is left and right adjoint to $H \mapsto \Res_{L/K} H$. In particular there are adjunction homomorphisms $\iota_{L/K}: G \to \Res_{L/K} G_L$ and $N_{L/K}:\Res_{L/K} G_L \to G$. When $L/K$ is a Galois extension, the Galois group $\Gal(L/K)$ acts canonically on $\Res_{L/K} G_L$. The following simple result, whose proof will be left to the reader, will be used later.

\begin{lemma}
\label{lemma:norm}
Let $L/K$ be a cyclic Galois extension of degree $n$, $\sigma$ a generator of $\Gal(L/K)$ and let $G$ be a commutative algebraic group over $K$. Let $G'$ be the kernel of $N_{L/K}:\Res_{L/K} G_L \to G$ so that $G'_L \cong G_L^{n-1}$. Then the map
\[
\Res_{L/K} (G_L)^{n-1} \cong\Res_{L/K}G'_L\stackrel{N_{L/K}}{\lra} G' \hookrightarrow \Res_{L/K} G_L 
\]
is given on the $i$-th summand by $1-\sigma^i$.
\end{lemma}

We denote by $\cC_K$ the category of finite reduced $K$-schemes. Thus each object of $\cC_K$ is isomorphic to $\Spec( E_1 \times \ldots \times E_r)$ where $E_1, \ldots, E_r/K$ are finite field extensions. A commutative algebraic group $G$ over $K$ defines a {\it Mackey-functor}, i.e.\ a co- and contravariant functor $G: \cC_K \to \ab$ satisfying (i), (ii) below. If $f:X \to Y$ is a morphism we denote by $f_*: G(X)\to G(Y)$ and $f^*: G(Y)\to G(X)$ the homomorphisms induced by co- and contravariant functoriality respectively. 
\begin{itemize}
\item[(i)] If $X = X_1 \coprod X_2 \in \Obj(\cC_K)$ then $G(X) = G(X_1) \oplus G(X_2)$.
\item[(ii)]
If
\[
\begin{CD}
X' @>> f' > Y'\\
@VV g' V @VV g V\\
X @>> f > Y
\end{CD}
\]
is a cartesian square in $\cC_K$ then $g^* f_* = (f')_* (g')^*$.
\end{itemize}
If $K \subseteq E_1 \subseteq E_2$ are two finite field extensions and $f: \Spec E_2 \to \Spec E_1$ the corresponding map in $\cC_K$ then $f^*$ (resp.\ $f_*$) is given by $\iota_{E_2/E_1}: G(E_1) \to G(E_2)$ (resp.\ $N_{E_2/E_1}:  G(E_2)\to G(E_1)$).

\paragraph{Local symbols} We recall also the notion of a {\it local symbol} (\cite{serre} and \cite{somekawa}) for $G$. Let $X\to \Spec K$ is a proper non-singular algebraic curve (note that we do not assume that $X$ is connected). Let $K(X)$ denote the ring of rational functions on $X$ and $|X|$ the set of closed points of $X$. For $P \in |X|$ we denote by $K_P$ the quotient field of the completion $\wcO_{X,P}$ of $\cO_{X,P}$, by $v_P: K_P\to \bZ\cup \{\infty\}$ the normalized valuation and by $K(P)$ the residue field of $K_P$. The local symbol at $P$ is a homomorphism $\partial_P: (K_P)^* \otimes G(K_P) \to G(K(P))$. It is characterized by the following properties:
\begin{itemize}
\item[(i)] If $f\in (K_P)^*$ and $g\in G(\wcO_{X,P})$ then $\partial_P(f\otimes g) =  v_P(f) g(P)$. Here $g(P)$ is the image of $g$ under the canonical map $G(\wcO_{X,P}) \to G(K(P))$.
\item[(ii)] For $f\in K(X)^*$ and $g\in G(K(X))$ we have $\sum_{P \in |X|} \, N_{K(P)/K}(\partial_P(f\otimes g)) = 0$.
\end{itemize}

\paragraph{Milnor $K$-groups attached to commutative algebraic groups} Let $G_1, \ldots, G_m$ be commutative algebraic groups over $K$. In \cite{somekawa} Somekawa has introduced the Milnor $K$-group $K(K;G_1,\ldots, G_m)$ (actually Somekawa considered only the case of semiabelian varieties though his construction works for arbitrary commutative algebraic groups). It is given as
\[
K(K;G_1,\ldots, G_m) = \left(\bigoplus_{X} \, G_1(X) \otimes \ldots \otimes G_m(X)\right)/R
\] 
where $X$ runs through all objects of $\cC_K$ and where the subgroup $R$ is generated by the following elements:

\noi (R1) If $f:X \to Y$ is a morphism in $\cC_K$ and if $x_{i_0}\in G_{i_0}(Y)$ for some $i_0$ and $x_i\in G_i(X)$ for $i\ne i_0$, then
\[
x_1\otimes \ldots\otimes f_*(x_{i_0}) \otimes\ldots \otimes x_m - f^*(x_1)\otimes \ldots\otimes x_{i_0} \otimes\ldots \otimes f^*(x_m)\in R.
\]

\noi (R2) Let $X\to \Spec K$ be a proper non-singular curve, $f\in K(X)^*$ and $g_i\in G_i(K(X))$. Assume that for each $P\in |X|$ there exists $i(P)$ such that $g_i\in G_i(\wcO_{X,P})$ for all $i\ne i(P)$. Then 
\[
\sum_{P \in |X|} \, g_1(P) \otimes \ldots \otimes \partial_P(f\otimes g_{i(P)}) \otimes \ldots \otimes g_m(P)\in R.
\]

For $X\in \cC_K$ and $x_i\in G_i(X)$ for $i=1, \ldots , m$ we write $\{x_1, \ldots, x_m\}_{X/K}$ for the image of $x_1\otimes \ldots \otimes x_m$ in $K(K;G_1,\ldots, G_m)$ (elements of this form will be referred to as {\it symbols}).

A sequence of algebraic groups $G' \to G \to G''$ over $K$ will be called {\it Zariski exact} if $G'(E) \to G(E) \to G''(E)$ is exact for every extension $E/K$. The proof of the following result is straightforward; hence will be omitted.  

\begin{lemma}
\label{lemma:rightexact}
Let $m$ be a positive integer and let $i\in \{1, \ldots, m\}$. Let $G_1, \ldots, G_m$ be commutative algebraic groups over $K$ and let $G_i' \to G_i \to G_i''\to 1$ be a Zariski exact sequence of commutative algebraic groups over $K$. Then the sequence
\[
K(K;G_1,\ldots, G_i',\ldots )\to K(K;G_1,\ldots, G_i,\ldots )\to K(K;G_1,\ldots, G_i'',\ldots)\to 0
\]
is exact as well.
\end{lemma}

\paragraph{The norm map}

Let $G_1, \ldots, G_m$ be commutative algebraic groups over $K$ and let $L/K$ be a finite extension. Set $K(L;G_1,\ldots, G_m) \colon = K(L;(G_1)_{L},\break \ldots, (G_m)_{L})$. Then we have the norm map \cite{somekawa}
\begin{equation}
\label{eqn:norm}
N_{L/K}: K(L;G_1,\ldots, G_m)  \lra K(K;G_1,\ldots, G_m)
\end{equation}
defined on symbols by $N_{L/K}(\{x_1, \ldots, x_m\}_{X/L}) = \{x_1, \ldots, x_m\}_{X/K}$ for any $X\in \cC_L$ and $x_i \in G_i(X) ~(i=1, \ldots, m).$
We give another interpretation of (\ref{eqn:norm}) below when $L/K$ is separable. It is based on the following result.

\begin{lemma}
\label{lemma:weilres} 
Let $L/K$ be a finite separable extension and let $i, m$ be positive integers with $i\le m$. Let $G_1, \ldots, G_{i-1},  G_{i+1}, \ldots, G_m$ be commutative algebraic groups over $K$ and let $G_i$ be a commutative algebraic group over $L$. Then, we have an isomorphism
\[
K(K;G_1,\ldots, \Res_{L/K} G_i,\ldots, G_m ) \cong K(L;(G_1)_{L},\ldots, G_i , \ldots, (G_m)_{L}).
\]
\end{lemma}

{\em Proof.} To simplify the notation we assume that $i=m$. We denote by $\pi^{-1}: \cC_K \to \cC_L$ and $\pi: \cC_L \to \cC_K$ the functors  
\[
\pi^{-1}(X \to \Spec K) \colon = (X\otimes_K L \to \Spec L), 
\]
\[
\pi(Y \to \Spec L) \colon = ( Y \to \Spec L \to \Spec K).
\]
$\pi$ is left adjoint to $\pi^{-1}$. For $X\in \cC_K$ and $Y\in \cC_L$ let
\[
p_X: X\otimes_K L \lra X, \quad \io_Y: Y \lra Y\otimes_K L.
\]
be the adjunction morphisms. We define homomorphisms
\[
\phi: K(K;G_1,\ldots, G_{m-1}, \Res_{L/K} G_m) \lra K(K;(G_1)_{L},\ldots, (G_{m-1})_{L}, G_m),
\]
\[
\psi: K(L;(G_1)_{L},\ldots, (G_{m-1})_{L}, G_m) \lra  K(K;G_1,\ldots, G_{m-1}, \Res_{L/K} G_m). 
\]
as follows. For $X\in \cC_K$, $x_1\in G_1(X), \ldots,  x_{m-1}\in G_{m-1}(X)$ and $x_m\in G_m(X\otimes_K L)$ we put
\[
\phi(\{x_1, \ldots, x_m\}_{X/K}) = \{p^*(x_1), \ldots, p^*(x_{m-1}), x_m\}_{(X\otimes_K L)/L}.
\]
Conversely, for $Y\in \cC_L$ and $y_1\in G_1(Y), \ldots, y_m \in G_m(Y)$ let
\[
\psi(\{y_1, \ldots,  y_{m-1}, y_m\}_{Y/L} = \{y_1, \ldots, y_{m-1}, \io_*(y_m)\}_{Y/K}.
\]
One can easily verify that these maps are well-defined and mutually inverse to each other.
\enddemo

Let $G_1, \ldots, G_m$ be commutative algebraic groups over $K$ and let $L/K$ be a finite separable extension. Take any $i\in \{1,\ldots, m\}$. The map $N_{L/K}: \Res_{L/K}(G_{i})_{L}\to G_i$ induces a map $K(K;G_1,\ldots, \Res_{L/K} (G_{i})_{L},\ldots, G_m) \lra K(K;G_1,\ldots, G_m)$, and the composition of it with the isomorphism $\psi$ above coincides with the norm map \eqref{eqn:norm}. When $L/K$ is a Galois extension, the action of $\Gal(L/K)$ on $\Res_{L/K} (G_{i})_{L}$ induces its action on
\[
K(L;G_1,\ldots, G_m)\cong K(K;G_1,\ldots, \Res_{L/K} (G_{i})_{L},\ldots, G_m)
\]
and we have $N_{L/K} \circ \sigma = N_{L/K}$ for all $\sigma\in \Gal(L/K)$. This action does not depend on the choice of $i$.

\begin{lemma}
\label{lemma:hilbert90} 
Let $L/K$ be a cyclic Galois extension and let $\sigma\in \Gal(L/K)$ be a generator. Suppose that for two different $i\in \{1, \ldots, m\}$ the sequence 
\begin{equation}
\label{eqn:normsurj}
\Res_{L/K}(G_{i})_{K}\stackrel{N_{L/K}}{\lra} G_i\lra 1
\end{equation} 
is Zariski exact. Then the sequence of abelian groups 
\[
K(L;G_1,\ldots, G_m)\stackrel{1-\sigma}{\lra} K(L;G_1,\ldots, G_m)\stackrel{N_{L/K}}{\lra} K(K;G_1,\ldots, G_m)\lra 0
\]
is exact. 
\end{lemma}

{\em Proof.} Suppose that (\ref{eqn:normsurj}) is exact for $i=m-1, m$. Let $G_m' \colon = \Ker(N_{L/K}:\Res_{L/K}(G_m)_L\to G_m)$. By Lemmas \ref{lemma:rightexact} and \ref{lemma:weilres} there are exact sequences
\begin{align}
\label{align:exakt1}
& K(K;G_1,\ldots, G_m')\to K(L;G_1,\ldots, G_m)\stackrel{N_{L/K}}{\lra} K(K;G_1,\ldots, G_m)\to 0\\
\label{align:exakt2}
& K(L;G_1,\ldots, G_{m-1}, G_m')\stackrel{N_{L/K}}{\lra} K(K;G_1,\ldots, G_{m-1}, G_m')\lra 0.
\end{align} 
Since $(G_m')_L\cong (G_m)_L^{n-1}$ ($n\colon =[L:K]$) we can replace the first group of (\ref{align:exakt2}) by $K(L;G_1,\ldots, G_m)^{n-1}$. By Lemma \ref{lemma:norm} the composite
\[
K(L;G_1,\ldots, G_m)^{n-1} \to K(K;G_1,\ldots, G_{m-1}, G_m')\to K(L;G_1,\ldots, G_m)
\]
is given on the $i$-th summand by $1-\sigma^i$. The assertion follows. \enddemo

\paragraph{Galois symbol}

Let $G_1, \ldots, G_m$ be connected commutative algebraic groups over $K$, and let $n$ be an integer prime to the characteristic of $K$. 
For any finite extension $L/K$, we have a homomorphism \cite{somekawa}
\begin{equation}
\label{eqn:galoissymbol}
h_L: K(L; G_1, \ldots, G_m)/n 
\lra H^m(L,G_1[n]\otimes\ldots\otimes G_m[n])
\end{equation}
called the {\it Galois symbol}.
This is characterized by the following properties.
\begin{itemize}
\item[(i)] 
If $x_i \in G_i(L)$ for $i=1, \ldots, m$,
then
$h_L( \{x_1, \ldots, x_m \}_{L/L} ) = (x_1) \cup \ldots \cup (x_m)$.
Here we write by $(x_i)$ for the image of
$x_i$ in $H^1(L, G_i[n])$ by the connecting homomorphism
associated to the exact sequence 
$1 \to G_i[n] \to G_i \overset{n}{\to} G_i \to 1.$
\item[(ii)] 
If $M/L/K$ is a tower of finite extensions and if $M/L$ is
separable (resp. purely inseparable), then the diagram 
\[
\begin{CD}
K(M; G_1, \ldots, G_m)/n @>> h_M > 
H^m(M,G_1[n]\otimes\ldots\otimes G_m[n])\\
@VV N_{M/L} V 
@VVV\\
K(L; G_1, \ldots, G_m)/n @>> h_L > 
H^m(L, G_1[n]\otimes\ldots\otimes G_m[n])\\
\end{CD}
\]
is commutative, where the right vertical map is the corestriction
(resp. the multiplication by $[M:L]$
under the identification
$H^m(M,G_1[n]\otimes\ldots\otimes G_m[n]) \cong
H^m(L, G_1[n]\otimes\ldots\otimes G_m[n])$ ).
\end{itemize}
Property (i) implies in particular that (\ref{eqn:galoissymbol}) coincides with the usual Galois symbol (\ref{eqn:galoisclass}) in the case $G_1= \ldots = G_m = \bG_m$. In \cite{somekawa} Remark 1.7, Somekawa conjectured that the Galois symbol associated to semiabelian varieties should be injective.

\paragraph{Galois cohomology of cyclic extensions}

Let $L/K$ be a cyclic Galois extension of degree $n$ and let $\sigma$ be a generator of $G \colon = \Gal(L/K)$. For a discrete $G_K$-module $M$, tensoring the short exact sequence of $G$-modules
\begin{equation}
\label{eqn:char}
0\lra \bZ \lra \bZ[G] \stackrel{1-\sigma}{\lra} \bZ[G] \lra \bZ\lra 0
\end{equation}
with $M$ yields a distinguished triangle 
\begin{equation}
\label{eqn:triangle}
M[1] \stackrel{\alpha}{\lra} C^{\punkt}(M) \stackrel{\beta}{\lra} M \stackrel{\gamma}{\lra} M[2]
\end{equation}
in the derived category $D(G_K)$. Here we denote by $C^{\punkt}(M)$ the complex
\[
\Res_{L/K} M \stackrel{1-\sigma}{\lra}\Res_{L/K} M
\]
concentrated in degree $-1$ and $0$. The spectral sequence
\[
E^{p,q}_1 = H^q(K, C^p(M))\Longrightarrow E^{p+q} = H^{p+q}(K,C^{\punkt}(M))
\]
induces short exact sequences
\begin{equation}
\label{eqn:edge}
0 \to H^q(L, M)_G \to H^q(K,C^{\punkt}(M)) \to H^{q+1}(L, M)^G\to 0.
\end{equation}
It is easy to see that the composite
\[
H^{q+1}(K, M)\stackrel{\alpha}{\lra}H^q(K,C^{\punkt}(M))\to H^{q+1}(L, M)^G
\]
is the restriction and 
\[
H^q(L, M)_G \to H^q(K,C^{\punkt}(M)) \stackrel{\beta}{\lra}H^q(K,M)
\]
is induced by the corestriction. In particular we have $\gamma(H^q(K,M)) \subseteq \break \Ker(\res: H^{q+2}(K, M) \to H^{q+2}(L, M))$ hence 
\begin{equation}
\label{eqn:killn}
n \gamma(H^q(K,M)) = 0.
\end{equation}

For an integer $m$ prime to $\charac K$ and $r\in \bN$ we write $\bZ/m\bZ(r) \colon = \mu_m^{\otimes^r}$ and 
\[
H^3(L/K,\bZ/m\bZ(2))\colon =  \Ker(H^3(K, \bZ/m\bZ(2))\stackrel{\res}{\lra}H^3(L, \bZ/m\bZ(2))).
\]
By restricting $\alpha: H^3(K, \bZ/m\bZ(2)) \to H^2(K,C^{\punkt}(\bZ/m\bZ(2)))$ to the subgroup $H^3(L/K,\bZ/m\bZ(2))$ and composing it with the inverse of the first map in (\ref{eqn:edge}) we obtain a map
\begin{equation}
\label{eqn:keymap}
 H^3(L/K,\bZ/m\bZ(2)) \to \Ker(H^2(L, \bZ/m\bZ(2))_G\stackrel{\cor}{\to} H^2(K, \bZ/m\bZ(2))).
\end{equation}

\begin{lemma}
\label{lemma:keylemma}
Assume that $n$ is prime to $\charac K$ and $\mu_{n^2}(\overline{K}) \subset K$. Then the homomorphism (\ref{eqn:keymap}) is injective for $m=n$.
\end{lemma}

{\em Proof.} It is enough to show that $\gamma: H^1(K, \bZ/n\bZ(2)) \to H^3(K, \bZ/n\bZ(2))$ is zero. Consider the commutative diagram 
\[
\begin{CD}
H^1(K, \bZ/n\bZ(2)) @>>>H^1(K, \bZ/n^2\bZ(2))\\
@VV \gamma V@VV \gamma V\\
H^3(K, \bZ/n\bZ(2)) @>>>H^3(K, \bZ/n^2\bZ(2))
\end{CD}
\]
induced by the canonical injection $\bZ/n\bZ(2)\to \bZ/n^2\bZ(2)$. The assumption $\mu_{n^2}(\overline{K}) \subset K$ implies that the upper horizontal map can be identified with
\[
K^*/(K^*)^n \lra K^*/(K^*)^{n^2}, x (K^*)^n \mapsto x^n (K^*)^{n^2}.
\]
In particular the image is contained in $n H^1(K, \bZ/n^2\bZ(2))$. By (\ref{eqn:killn}) it is mapped under $\gamma$ to $n\gamma(H^1(K, \bZ/n^2\bZ(2)))= 0$. On the other hand it is a simple consequence of the Merkurjev-Suslin theorem \cite{ms} that the lower horizontal map is injective. Hence $\gamma(H^1(K, \bZ/n(2))) = 0$. \enddemo

\paragraph{The counterexample}

Let $L/K$ be as in the last section and let $T \colon = \Ker(N_{L/K}:\Res_{L/K} \bG_m \to \bG_m)$. We make the following assumptions
\begin{align}
\label{align:torus1}
& \mbox{$n$ is prime to $\charac K$ and $\mu_{n^2}(\overline{K}) \subset K$,}\hspace{3cm}\\
\label{align:torus2}
& \mbox{$H^3(L/K,\bZ/n\bZ(2)))\ne 0$.}
\end{align}

\begin{proposition}
\label{proposition:counterex} 
The Galois symbol $K(K; T,T)/n \to H^2(K, T[n]\otimes T[n])$ is not injective.
\end{proposition}

{\em Proof.} Let $\sigma$ be a generator of $G \colon = \Gal(L/K)$. The exact sequence  
\[
1\lra \bG_m \lra\Res_{L/K} \bG_m \stackrel{1-\sigma}{\lra}\Res_{L/K} \bG_m\stackrel{N_{L/K}}{\lra} \bG_m \lra 1
\]
yields two short exact sequences 
\begin{align}
\label{align:torus3}
& 1\lra \bG_m \lra \Res_{L/K} \bG_m \lra T \lra 1,\hspace{2cm} \\ 
\label{align:torus4} 
& 1 \lra T \lra \Res_{L/K} \bG_m\lra \bG_m\lra 1.
\end{align}

Correspondingly, (\ref{eqn:char}) induces two short exact sequences 
\begin{equation}
\label{eqn:cochar}
0 \to \bZ \to \bZ[G] \to X \to 0, \qquad 0 \to X \to\bZ[G] \to\bZ \to 0
\end{equation}
where $X$ denotes the cocharacter group of $T$. Note that the sequence (\ref{align:torus3}) is Zariski exact by Hilbert 90. Since the map $\Res_{L/K} \bG_m \to T$ factors through $\Res_{L/K} \bG_m \to \Res_{L/K} T \to T$ the sequence $\Res_{L/K} T \to T \to 1$ is Zariski exact as well. By Lemma \ref{lemma:hilbert90} the upper horizontal map in the diagram
\[
\begin{CD}
(K(L; T, T)/n)_G  @> N_{L/K} >> K(K; T,T)/n\\
@VVV@VVV\\
H^2(L, T[n]\otimes T[n])_G @> \cor >>H^2(K, T[n]\otimes T[n])
\end{CD}
\]
is an isomorphism. The vertical maps are Galois symbols. Since $T_L$ is a split torus the left vertical map is an isomorphism by the Merkurjev-Suslin theorem \cite{ms}. Thus to finish the proof it remains to show that the lower vertical arrow is not injective. Note that $T[n] \cong \bZ/n\bZ(1)\otimes X$. Hence the assertion follows from Lemma \ref{lemma:keylemma} and Lemma \ref{lemma:directsum} below.

\begin{lemma}
\label{lemma:directsum}
There exists homomorphisms of $G$-modules $e: \bZ \to X\otimes_{\bZ} X$ and $f: X\otimes_{\bZ} X\to \bZ$ such that $f\circ e: \bZ\to \bZ$ is multiplication by $n-1$.
\end{lemma}

{\em Proof.} For a $G$-module $M$ we write $M^{\vee}$ for the $G$-module $\Hom(M, \bZ)$. Let $(\,\,,\,\,):\bZ[G] \otimes_{\bZ} \bZ[G] \lra \bZ$
be the symmetric pairing given by 
\begin{equation}
\label{eqn:pair}
(g,g') \quad =  \quad  \left\{ \begin{array}{ll} 
                       1 & \mbox{if $g=g'$,}\\
                       0 & \mbox{if $g\ne g'$.}\\
\end{array} \right.
\end{equation}
It yields an isomorphism $\bZ[G]\to \bZ[G]^{\vee}$. For a submodule $M \subseteq \bZ[G]$ let 
\[
M^{\perp} = \{x\in \bZ[G]\mid \, (x,m) = 0 \,\,\,\forall\, m\in M\}.
\]
Then we have $X^{\perp} = \bZ S$ and $(\bZ S)^{\perp} = X$ where $S = \sum_{i=0}^{n-1} \, \sigma^i$. Thus (\ref{eqn:pair}) yields an isomorphism
$X \cong (\bZ[G]/\bZ S)^{\vee}$. By (\ref{eqn:cochar}) we have $\bZ[G]/\bZ S\cong X$, hence 
\[
X\otimes_{\bZ} X\quad \cong\quad X\otimes_{\bZ} X^{\vee} \quad \cong\quad \Hom(X,X)
\]
Thus it suffices to prove the assertion for $\Hom(X,X)$. Obviously, for the two maps $e:\bZ \to \Hom(X,X), m\mapsto m\id_X$ and $f: \Hom(X,X) \to \bZ, \tau \mapsto \Tr(\tau)$ we have $f\circ e = \rank(X) = n-1$. \enddemo

\begin{remark}
\label{remark:nlocal} \rm It is easy to construct examples where the assumptions (\ref{align:torus1}) and (\ref{align:torus2}) above are satisfied. For instance if $K$ is a $2$-local field satisfying property (\ref{align:torus1}) and $L/K$ is any cyclic extension of degree $n$ then (\ref{align:torus2}) holds by \cite{kato}.
\end{remark}

\section{Counterexample to a conjecture of Beilinson}
\label{section:beilinson}
 
We first introduce some notation and recall a few facts from \cite{voevodsky} and \cite{mvw}. Let $K$ be a field of characteristic zero. Let $Cor_K$ denote the additive category of finite correspondences (\cite{mvw}, 1.1). The objects of $Cor_K$ are smooth separated $K$-schemes of finite type and for $X,Y\in \Obj(Cor_K)$ the group of morphisms $Cor_K(X,Y)$ is the free abelian group generated by integral closed subschemes $W$ of $X\times Y$ which are finite and surjective over $X$. Let $\Dnis{K}$ (resp.\ $\Det{K}$) denote the derived category of complexes of Nisnevich (resp.\ {\'e}tale) sheaves with transfer bounded from above. 

The category of effective motivic complexes $\DMn(K)$ (resp.\ {\'e}tale effective motivic complexes $\DMe(K)$) is the full subcategory of $D^{-}(Shv_{\Nis}\break (Cor_K))$ (resp.\ $\Det{K}$) which consists of complexes $C^\star$ with homotopy invariant cohomology sheaves $H^i(C^\star)$ for all $i$ (see \cite{voevodsky}, 3.1 or \cite{mvw}, 14.1, resp.\ 9.2). $\DMn(K)$ and $\DMe(K)$ are triangulated tensor categories. They are equipped with the t-structure induced from the standard t-structure on $\Dnis{K}$ (resp.\ $\Det{K}$). There is a covariant functor $M: Cor_K\to  \DMn(K), X \mapsto M(X)$ and we have $M(X\times Y) = M(X)\otimes M(Y)$. There is also the "change of topology" functor $\alpha^*: \DMn(K)\to \DMe(K)$. It is a tensor functor which admits a right adjoint $R\alpha_*: \DMe(K)\to \DMn(K)$. 

Beilinson \cite{beilinson} has proposed the following generalization of the Milnor-Bloch-Kato conjecture: For any smooth affine $K$-scheme $X$ the adjunction morphism $M(X)\to R\alpha_*\alpha^* M(X)$ induces an isomorphism on cohomology in degrees $\le 0$, i.e.\ the map 
\begin{equation}
\label{eqn:beilinsonconj}
a_X: M(X)\lra t_{\le 0} R\alpha_*\alpha^* M(X)
\end{equation}
is an isomorphism in $\DMn(K)$. 

If $X = (\bG_m)^d= \bG_m\times \ldots \times \bG_m$ ($d$-fold product of $\bG_m$) we have $M(X) \cong (\bZ \oplus \bZ(1)[1])^d$. Thus $a_X$ is an isomorphism if and only if 
\begin{equation}
\label{eqn:mbkconj}
\bZ(n) \lra t_{\le n} R\alpha_*\alpha^* \bZ(n)
\end{equation}
is an isomorphism for all $n\le d$. It is known (compare \cite{sv}) that the Milnor-Bloch-Kato conjecture is equivalent to the assertion that (\ref{eqn:mbkconj}) is an isomorphism for all $n\ge 0$.

Let $L/K$ be a separable quadratic extension and let $T \colon = \Ker(N_{L/K}:\Res_{L/K} \bG_m \to \bG_m)$. We shall show that (\ref{eqn:beilinsonconj}) is in general not an isomorphism for $X = T^n$ for $n\ge 2$.
By (\cite{hk}, 7.3) there exists a canonical decomposition $M(T) = \bZ \oplus \bZ(L/K, 1)[1]$
where $\bZ(L/K, 1)$ is the cone of the morphism $\bZ(1)\to \Res_{L/K} \bZ(1)$. 

\begin{remarks}
\rm (a) Here is a more explicit description of the motive $\bZ(L/K, 1)$. The torus $T$ defines a homotopy invariant {\'e}tale (hence Nisnevich) sheaf with transfer and therefore an element of $\DMn(K)$. We have 
\[
\bZ(L/K, 1) \quad \cong \quad T[-1].
\]
This can be deduced from the corresponding statement for $\bG_m$ (\cite{mvw}, 4.1) and the exactness of (\ref{align:torus3}) (as a sequence in $Shv_{\Nis}(Cor_K)$).

\noi (b)\footnote{This remark has been communicated to us by B.\ Kahn.} Let $A_1, \ldots , A_n$ be semi-abelian varieties over $K$. It should be possible to identify the generalized Milnor $K$-group $K(K;A_1,\ldots, A_n)$ with a $\Hom$-group in $\DMn(K)$. For that we view $A_1, \ldots , A_n$ again as elements in $Shv_{\Nis}(Cor_K)$. Then we expect that
\[
K(K;A_1,\ldots, A_n) \cong \Hom_{\DMn(K)}(\bZ, A_1\otimes \ldots \otimes A_n).
\]
If $A_1= \ldots = A_n = \bG_m$ this is proved in (\cite{mvw}, lecture 5) and it is likely that the proof given there can be adapted to the case of arbitrary semi-abelian varieties. 
\end{remarks}

For $p,q \ge 0$ and $n= p+q$ we define
\[
\bZ(L/K, p, q) \colon = \bZ(L/K, 1)^{\otimes^p}\otimes \bZ(q)
\]
and denote by $C(p,q)$ the cone of $\bZ(L/K, p, q) \lra t_{\le n} R\alpha_*\alpha^* \bZ(L/K, p, q)$.
Note that $\bZ(L/K, p, q)[n]$ is a direct summand of $M(T^p\times (\bG_m)^q)$. We also put $C(n) \colon = C(0,n)$. We have
\begin{equation}
\label{eqn:mbk2}
C(n)\cong  (t_{\ge n+1} R\alpha_* \bQ/\bZ(n))[-1]
\end{equation}
This follows from the Milnor-Bloch-Kato conjecture (in fact for our purpose we need (\ref{eqn:mbk2}) only after localization at the prime $2$ where it follows from the Milnor conjecture \cite{voevodsky2}).

Tensoring $\bZ(1)\to \Res_{L/K} \bZ(1)\to \bZ(L/K,1) \to \bZ(1)[1]$ with $\bZ(L/K, p-1, q)$ (for $p\ge 1, q\ge 0$) yields a distinguished triangle
\[
\bZ(L/K, p-1, q+1) \to \Res_{L/K} \bZ(n) \to \bZ(L/K, p, q)\to \bZ(L/K, p-1, q+1)[1]
\]
hence also a triangle
\begin{equation}
\label{eqn:basictri}
C(p-1, q+1)\to \Res_{L/K} C(n) \to C(p,q)\to C(p-1, q+1)[1].
\end{equation}
The following Lemma follows easily by induction on $q$ using (\ref{eqn:mbk2}) and (\ref{eqn:basictri}).

\begin{lemma}
\label{lemma:counterexample}
Let $p\ge 1,q\ge 0$ and $n= p+q$. Then we have $H^{k}(C(p,q)) = 0$ for $k<q+2$ and
\[
H^{q+2}(C(p,q))(K) \cong H^{n+1}(L/K, \bQ/\bZ(n))
\]
where $H^{n+1}(L/K, \bQ/\bZ(n))\colon = \Ker(H^{n+1}(K, \bQ/\bZ(n))\stackrel{\res}{\lra}H^{n+1}(L, \bQ/\bZ(n)))$.
\end{lemma}

Since $[L:K]=2$ we have 
\begin{align*}
H^{n+1}(L/K, \bQ/\bZ(n)) & \cong H^{n+1}(L/K, \bQ_2/\bZ_2(n)) \cong H^{n+1}(L/K, \bZ/2\bZ(n))\\
& \cong H^{n+1}(L/K, \bZ/2\bZ)
\end{align*}
(the second isomorphism is a consequence of the Milnor conjecture).
Now the following Proposition follows by 
applying Lemma \ref{lemma:counterexample} for $(p,q) = (2,0)$ and $(n, 0)$.

\begin{proposition}
\label{proposition:beilcounterexample} (a) There exists a short exact sequence
\[
0\lra H^0(M(T\times T))(K)\lra R^0\alpha_*\alpha^* M(T\times T)(K)\lra H^3(L/K, \bZ/2\bZ)\lra 0
\]
In particular if $H^3(L/K, \bZ/2\bZ)\ne 0$ then (\ref{eqn:beilinsonconj}) is not an isomorphism for $X = T\times T$.

\noi (b) More generally let $n$ be an integer $\ge 2$ and assume that $H^{n+1}(L/K, \bZ/2\bZ)\ne 0$. Then the map (\ref{eqn:beilinsonconj}) is not an isomorphism for $X = T^n$. More precisely either the map
\[
H^{2-n}(M(X))\lra R^{2-n}\alpha_*\alpha^* M(X)
\]
is not surjective or
\[
H^{3-n}(M(X))\lra R^{3-n}\alpha_*\alpha^* M(X)
\]
is not injective. 
\end{proposition}

An $n$-local field $K$ of characteristic $0$ provides an example where the above assumption holds. In fact by \cite{kato} we have $H^{n+1}(L/K, \bZ/2\bZ)\cong \bZ/2\bZ$ for such fields.

\bigskip
\begin{tabbing} 
 \hspace{1cm}\= Michael Spie{\ss}\hspace{4cm} \= Takao Yamazaki\\
 \>Fakult{\"a}t f{\"u}r Mathematik \>Mathematical Institute\\
 \>Universit{\"a}t Bielefeld \> Tohoku University\\
\>Postfach 100131 \> Aoba\\
\>D-33501 Bielefeld, Germany \> Sendai 980-8578, Japan\\
\> mspiess@math.uni-bielefeld.de \> ytakao@math.tohoku.ac.jp
\end{tabbing}

\end{document}